\documentclass[12pt]{article}%
\usepackage{amsfonts}
\usepackage{amssymb}
\usepackage{amsmath}%
\usepackage{graphicx}
\providecommand{\U}[1]{\protect\rule{.1in}{.1in}}
\newtheorem{theorem}{Theorem}

\newtheorem{corollary}[theorem]{Corollary}

\begin{document}

\title{Choosing elements from finite fields}
\author{Michael Vaughan-Lee}
\date{November 2012}
\maketitle

\section{Introduction}

Graham Higman wrote two immensely important and influential papers on
enumerating $p$-groups in the late 1950s. The papers were entitled
\emph{Enumerating }$p$\emph{-groups} I and II, and were published in the
Proceedings of the London Mathematical Society in 1960 (see \cite{higman60}
and \cite{higman60b}). In these two papers Higman proved that for any given
$n$, the function $f(p^{n})$ enumerating the number of $p$-groups of order
$p^{n}$ is bounded by a polynomial in $p$, and he formulated his famous PORC
conjecture concerning the form of the function $f(p^{n})$. He conjectured that
for each $n$ there is an integer $N$ (depending on $n$) such that for $p $ in
a fixed residue class modulo $N$ the function $f(p^{n})$ is a polynomial in
$p$. For example, {for }$p\geq5$ the number of groups of order $p^{6}$ is{%
\[
3p^{2}+39p+344+24\gcd(p-1,3)+11\gcd(p-1,4)+2\gcd(p-1,5).
\]
(See \cite{newobvl}.) So for }$p\geq5$, $f(p^{6})$ is one of 8 polynomials in
$p$, with the choice of polynomial depending on the residue class of $p$
modulo 60. The number of groups of order $p^{6}$ is \textbf{P}olynomial
\textbf{O}n \textbf{R}esidue \textbf{C}lasses. As evidence in support of his
PORC conjecture Higman proved that, for any given $n$, the function
enumerating the number of $p$-class 2 groups of order $p^{n}$ is a PORC
function of $p$. He obtained this result as a corollary to a very general
theorem about vector spaces acted on by the general linear group. As another
corollary to this general theorem, he also proved that for any given $n$ the
function enumerating the number of algebras of dimension $n$ over the field of
$q$ elements is a PORC function of $q$. A key step in Higman's proof of these
results is Theorem 2.2.2 from \cite{higman60b}.

\begin{theorem}
[Higman \cite{higman60b}]The number of ways of choosing a finite number of
elements from $\mathbb{F}_{q^{n}}$ subject to a finite number of monomial
equations and inequalities between them and their conjugates over
$\mathbb{F}_{q}$, considered as a function of $q$, is PORC.
\end{theorem}

The statement of this theorem probably requires some explanation! Here we are
choosing elements $x_{1},x_{2},\ldots,x_{k}$ (say) from the finite field
$\mathbb{F}_{q^{n}}$ (where $q$ is a prime power) subject to a finite set of
equations and non-equations of the form%
\[
x_{1}^{n_{1}}x_{2}^{n_{2}}\ldots x_{k}^{n_{k}}=1
\]
and%
\[
x_{1}^{n_{1}}x_{2}^{n_{2}}\ldots x_{k}^{n_{k}}\neq1,
\]
where $n_{1},n_{2},\ldots,n_{k}$ are integer polynomials in the Frobenius
automorphism $x\mapsto x^{q}$ of $\mathbb{F}_{q^{n}}$. Higman calls these
equations and non-equations \emph{monomial}. For example, suppose we want to
choose $x_{1},x_{2}\in\mathbb{F}_{q^{n}}$ such that $x_{1}$ is the root of an
irreducible quadratic over $\mathbb{F}_{q}$ and such that $x_{2}^{2}$ is the
product of the roots of this quadratic. Then we require that $x_{1}$ and
$x_{2}$ satisfy%
\begin{equation}
x_{1}^{q^{2}-1}=1,\;x_{1}^{q-1}\neq1,\;x_{1}^{q+1}x_{2}^{-2}=1.
\end{equation}
The equation $x_{1}^{q^{2}-1}=1$ guarantees that $x_{1}$ is the root of a
quadratic over $\mathbb{F}_{q}$, and the non-equation $x_{1}^{q-1}\neq1$
guarantees that $x_{1}\notin\mathbb{F}_{q}$ so that the quadratic is
irreducible. The other root of the quadratic is then $x_{1}^{q}$, so the last
equation guarantees that $x_{2}^{2}$ is the product of the roots. To make sure
that $x_{1},x_{2}\in\mathbb{F}_{q^{n}}$, we also require that $x_{1}^{q^{n}%
-1}=1$, $x_{2}^{q^{n}-1}=1$. Higman's theorem is that the function enumerating
the number of solutions to (1) in $\mathbb{F}_{q^{n}}$ is a PORC function of
$q$.

In this note we give very precise information about the exact form of the PORC
functions needed to enumerate the number of solutions to a set of monomial equations.

Higman's proof of his Theorem 2.2.2 involves five pages of homological
algebra. A shorter more elementary proof can be found in \cite{vlee12a}. The
proof in \cite{vlee12a} shows that one way to calculate the number of
solutions to a set of monomial equations is to write the equations as the rows
of a matrix. So we represent the equations%
\[
x_{1}^{q^{2}-1}=1,\;x_{1}^{q+1}x_{2}^{-2}=1,\;x_{1}^{q^{n}-1}=1,\;x_{2}%
^{q^{n}-1}=1
\]
by the matrix%
\[
\left[
\begin{array}
[c]{cc}%
q^{2}-1 & 0\\
q+1 & -2\\
q^{n}-1 & 0\\
0 & q^{n}-1
\end{array}
\right]  .
\]
For any given value of $q$ the matrix becomes an integer matrix, and it is
shown in \cite{vlee12a} that the number of solutions is the product of the
elementary divisors in the Smith normal form of this integer matrix. To obtain
the number of solutions to (1) in $\mathbb{F}_{q^{n}}$ we subtract the number
of solutions to the equations%
\[
x_{1}^{q-1}=1,\;x_{1}^{q+1}x_{2}^{-2}=1,\;x_{1}^{q^{n}-1}=1,\;x_{2}^{q^{n}%
-1}=1.
\]
The number of solutions to these equations is just the product of the
elementary divisors in the Smith normal form of the matrix%
\[
\left[
\begin{array}
[c]{cc}%
q-1 & 0\\
q+1 & -2\\
q^{n}-1 & 0\\
0 & q^{n}-1
\end{array}
\right]  .
\]
(In \cite{vlee12a} $q$ is assumed to be prime, but the proof is still valid
when $q$ is a prime power.)

Matrices used in this way to represent a set of monomial equations have
entries which are integer polynomials in $q$. The columns correspond to the
unknowns we are solving for, and since there will always be rows
$(q^{n}-1,0,0,\ldots,0)$, $(0,q^{n}-1,0,\ldots,0)$, \ldots, $(0,\ldots
,0,q^{n}-1)$ corresponding to the requirement that the unknowns are elements
in $\mathbb{F}_{q^{n}}$, it follows that the rank of one of these matrices is
the number of columns. So the product of the elementary divisors in the Smith
normal form of one of these matrices is the greatest common divisor of the
$k\times k$ minors, where $k$ is the number of columns. These $k\times k$
minors are integer polynomials in $q$, and it is proved in \cite{vlee12a} that
the greatest common divisor of a set of integer polynomials in $q$ is a PORC
function of $q$. There is some ambiguity about what \textquotedblleft the
greatest common divisor of a set of polynomials\textquotedblright\ means here.
Suppose we have some integer polynomials $f_{1}(q),f_{2}(q),\ldots,f_{s}(q)$.
For any given value of $q$ these polynomials evaluate to integers, and by
\textquotedblleft greatest common divisor of the polynomials\textquotedblright%
\ we actually mean \textquotedblleft greatest common divisor of the values of
the polynomials at $q$\textquotedblright. It is this integer valued function
of $q$ which we claim is PORC, and it turns out that we can be quite precise
about the form that this PORC function takes.

\begin{theorem}
The greatest common divisor of a set of integer polynomials in $q$ can be
expressed in the form $df$ where $f$ is an integer polynomial in $q$ and where%
\[
d=\alpha+\sum_{i=1}^{r}\alpha_{i}\gcd(q-n_{i},m_{i})
\]
for some rational numbers $\alpha,\alpha_{1},\alpha_{2},\ldots,\alpha_{r}$,
some integers $m_{1},m_{2},\ldots,m_{r}$ with $m_{i}>1$ for all $i$, and some
integers $n_{1},n_{2},\ldots,n_{r}$ with $0<n_{i}<m_{i}$ for all $i $.
\end{theorem}

\begin{corollary}
The number of ways of choosing a finite number of elements from $\mathbb{F}%
_{q^{n}}$ subject to a finite number of monomial equations and inequalities
between them and their conjugates over $\mathbb{F}_{q}$ can be expressed as a
linear combination of terms of the form $df\,$, where $f$ and $d$ are as
described in Theorem 2.
\end{corollary}

\section{Choosing field elements}

To make this note self contained, we give a proof here that we can use a
matrix to represent a set of monomial equations over a finite field, and that
the number of solutions to the equations is the product of the elementary
divisors in the Smith normal form of the matrix.

So suppose we have a set of monomial equations in unknowns $x_{1},x_{2}%
,\ldots,x_{k}$, and suppose that we want to find the number of solutions to
these equations in the field $\mathbb{F}_{q^{n}}$. We represent the equations
in a matrix $A$ with $k$ columns, with a row $(n_{1},n_{2},\ldots,n_{k})$ for
each monomial equation $x_{1}^{n_{1}}x_{2}^{n_{2}}\ldots x_{k}^{n_{k}}=1$. We
also add in rows%
\[
(q^{n}-1,0,0,\ldots,0),\,(0,q^{n}-1,0,\ldots,0),\,\ldots,\,(0,0,\ldots
,0,q^{n}-1)
\]
corresponding to the requirement that $x_{1},x_{2},\ldots,x_{k}\in
\mathbb{F}_{q^{n}}$. Note that the entries in the matrix $A$ are integer
polynomials in $q$. We now take a particular value for $q$ so that the matrix
becomes a matrix with integer entries.

Let $\omega$ be a primitive element in $\mathbb{F}_{q^{n}}$, and write
$x_{i}=\omega^{m_{i}}$ for $i=1,2,\ldots,k\,$, taking the exponents $m_{i}$ as
elements in $\mathbb{Z}_{q^{n}-1}$. Then a row $(\beta_{1},\beta_{2}%
,\ldots,\beta_{k})$ in the matrix $A$ corresponds to a relation $\beta
_{1}m_{1}+\beta_{2}m_{2}+\ldots+\beta_{k}m_{k}=0$ which we require the
exponents $m_{i}$ to satisfy. The matrix $A$ can be reduced to Smith normal
form over $\mathbb{Z}$ by elementary row and column operations. As we apply
these operations, the relations encoded in the matrix change. But we show that
at each step the number of solutions to the relations stays constant.

This is clear for elementary row operations, since an elementary row operation
replaces the relations by an equivalent set of relations. So we need to
consider the effect of elementary column operations. We can consider the
$k$-tuples $(m_{1},m_{2},\ldots,m_{k})$ as elements in the additive group
$G=\mathbb{Z}_{q^{n}-1}\times\mathbb{Z}_{q^{n}-1}\times\ldots\times
\mathbb{Z}_{q^{n}-1}$. Let $A$ be one of these relation matrices, and let $B$
be the matrix obtained from $A$ after applying an elementary column operation.
For each such operation we define an automorphism $\sigma$ of $G$ with the
property that $g\in G$ satisfies the relations given by the rows of $A$ if and
only if $g\sigma$ satisfies the relations given by the rows of $B $. This
shows that the number of elements in $G$ satisfying the relations given by $A$
is the same as the number of elements in $G$ satisfying the relations given by
$B$. If the elementary column operation swaps two columns of $A$ then we let
$\sigma$ be the automorphism which swaps the corresponding entries in
$(m_{1},m_{2},\ldots,m_{k})$, and if the elementary column operation
multiplies a column by $-1$ we let $\sigma$ be the automorphism which
multiplies the corresponding entry in $(m_{1},m_{2},\ldots,m_{k})$ by $-1$.
Finally, if the elementary column operation subtracts $\alpha$ times column
$j$ from column $i$, then we let $\sigma$ be the automorphism which leaves all
the entries in $(m_{1},m_{2},\ldots,m_{k})$ fixed except for the $j$-th entry,
which it replaces by $m_{j}+\alpha m_{i}$.

The argument above shows that the number of $g\in G$ satisfying the original
set of relations given by the rows of $A$ is the same as the number of $g\in
G$ satisfying the relations given by the Smith normal form $A$. If the
elementary divisors in the Smith normal form are $d_{1},d_{2},\ldots,d_{k}$,
then $(m_{1},m_{2},\ldots,m_{k})$ is a solution to these equations if and only
if%
\[
d_{1}m_{1}=d_{2}m_{2}=\ldots=d_{k}m_{k}=0.
\]
Provided we can show that $d_{i}|q^{n}-1$ for all $i$, this shows that the
number of solutions is $d_{1}d_{2}\ldots d_{k}$, as claimed.

If $A$ is one of these relation matrices with $k$ columns, then the rows of
$A$ are elements in the free $\mathbb{Z}$-module $F=\mathbb{Z}^{k}$. We let
$R(A)$ denote the $\mathbb{Z}$-submodule of $F$ generated by the rows of $A$.
Our claim that $d_{i}|q^{n}-1$ for all $i$ amounts to the claim that
$(q^{n}-1)F\leq R(S)$, where $S$ is the Smith normal form of our initial
relation matrix. The Smith normal form is obtained from the initial matrix by
a sequence of elementary row and column operations, and we show that
$(q^{n}-1)F\leq R(B)$ for all the matrices $B$ generated in this sequence.

Let $A$ be the starting matrix. Then it contains rows
\[
(q^{n}-1,0,0,\ldots,0),\,(0,q^{n}-1,0,\ldots,0),\,\ldots,\,(0,0,\ldots
,0,q^{n}-1),
\]
so it is clear that $(q^{n}-1)F\leq R(A)$. Suppose that at some intermediate
stage in the reduction of $A$ to Smith normal form we have two matrices $B$
and $C$, where $C$ is obtained from $B$ by an elementary row operation or an
elementary column operation. We assume by induction that $(q^{n}-1)F\leq R(B)
$, and we show that this implies that $(q^{n}-1)F\leq R(C)$. This is clear if
$C$ is obtained from $B$ by an elementary row operation, since then
$R(B)=R(C)$. So consider the case when $C$ is obtained from $B$ by an
elementary column operation. This column operation corresponds to an
automorphism $\sigma$ of $F$, and if $r$ is a row of $B$ then the
corresponding row of $C$ is $r\sigma$. So $R(C)=R(B)\sigma$, and the fact that
$(q^{n}-1)F$ is a characteristic submodule of $F$ implies that $(q^{n}-1)F\leq
R(C)$.

This completes the proof that the number of solutions to the relations given
by the rows of the matrix is equal to the product of the elementary divisors
in the Smith normal form. As mentioned in the introduction, the product of the
elementary divisors in the Smith normal form of an integer matrix with $k $
columns and rank $k$ is the greatest common divisor of the $k\times k$ minors.
In the situation we are concerned with, these minors are integer polynomials
in $q$. So the number of solutions to our monomial equations is the greatest
common divisor of a set of integer polynomials in $q$. More precisely, we have
a set of integer polynomials in $q$, and for any given value of $q$ the number
of solutions to our monomial equations is the greatest common divisor over
$\mathbb{Z}$ of the values of these polynomials at $q$.

\section{Proof of Theorem 2}

Let $f_{1}(q)$, $f_{2}(q)$, \ldots, $f_{s}(q)$ be a set of integer polynomials
in $q$. We want to compute the function whose value at $q$ is the greatest
common divisor of the integers $f_{1}(q)$, $f_{2}(q)$, \ldots, $f_{s}(q)$. In
this section there is no requirement that $q$ be a prime power, and to make
this clear we define a function $h:\mathbb{Z}\rightarrow\mathbb{Z}$ by setting%
\[
h(x)=\gcd(f_{1}(x),f_{2}(x),\ldots,f_{s}(x))\text{ for }x\in\mathbb{Z}.
\]
It is the function $h$ we want to compute. As mentioned above, there is some
ambiguity about what we mean by \textquotedblleft the greatest common divisor
of $f_{1}(x)$, $f_{2}(x)$, \ldots, $f_{s}(x)$\textquotedblright. We now
exploit this ambiguity, and treat $x$ as an indeterminate and treat $f_{1}%
(x)$, $f_{2}(x)$, \ldots, $f_{s}(x)$ as elements of the Euclidean domain
$\mathbb{Q}[x]$.

We can use the Euclidean algorithm to compute the greatest common divisor
$f(x)$ of $f_{1}(x)$, $f_{2}(x)$, \ldots, $f_{s}(x)$ in $\mathbb{Q}[x]$ and we
can take $f(x)$ to be a primitive polynomial in $\mathbb{Z}[x]$. We then
obtain polynomials $g_{1},g_{2},\ldots,g_{s}\in\mathbb{Q}[x]$ such that%
\[
f_{1}g_{1}+f_{2}g_{2}+\ldots+f_{s}g_{s}=f.
\]
Let $m$ be the least common multiple of the denominators of the coefficients
in $g_{1},g_{2},\ldots,g_{s}$. Then for any given value of $x$ in $\mathbb{Z}%
$, the greatest common divisor of the integers $f_{1}(x)$, $f_{2}(x)$, \ldots,
$f_{s}(x)$ is $df(x)$ for some $d$ dividing $m$. Furthermore, as a function of
$x$, the value of $d$ at $x$ depends only on the residue class of $x$ modulo
$m$. We show that we can express $d(x)$ in the form%
\[
\alpha+\sum_{i=1}^{r}\alpha_{i}\gcd(x-n_{i},m_{i})
\]
described in the statement of Theorem 2. Furthermore we show that we can take
the integers $m_{i}$ to be of the form $\frac{m}{d_{i}}$ for some square free
divisors $d_{i}$ of $m$ with $d_{i}<m$. This shows that $h(x)=d(x)f(x)$ has
the form described in Theorem 2.

If $m=1$ then $d=1$ for all $x$, and we are done. So suppose that $m>1$ and
let $S$ be the set of prime factors of $m$. For each subset $T\subset S$ let%
\[
d_{T}=%
{\displaystyle\prod\limits_{p\in T}}
p,
\]
and consider the function $k:\mathbb{Z}\rightarrow\mathbb{Z}$ defined by%
\[
k(x)=\sum_{T\subset S}(-1)^{|T|}\gcd(x,\frac{m}{d_{T}}).
\]
Clearly the value of $k$ at any given value of $x$ depends only on the residue
class of $x$ modulo $m$. First consider the case when $x=m$.%
\[
k(m)=\sum_{T\subset S}(-1)^{|T|}\frac{m}{d_{T}}=m%
{\displaystyle\prod\limits_{p\in S}}
(1-\frac{1}{p})\neq0.
\]
Next suppose that $1\leq x<m$. Then there is at least one $p\in S$ with the
property that the power of $p$ dividing $x$ is less than the power of $p$
dividing $m$. Pick one such $p$ and let $U=S\backslash\{p\}$. Then%
\[
k(x)=\sum_{T\subset U}(-1)^{|T|}\left(  \gcd(x,\frac{m}{d_{T}})-\gcd
(x,\frac{m}{pd_{T}})\right)  =0,
\]
since $\gcd(x,\frac{m}{d_{T}})=\gcd(x,\frac{m}{pd_{T}})$ for all $T\subset U$.
So if we let $c=k(m)$ then $\frac{1}{c}k(x)$ takes values $0,0,\ldots,0,1$ as
$x$ takes values $1,2,\ldots,m$ modulo $m$. It follows that if $0<a<m$ then
$\frac{1}{c}k(x-a)$ takes values $0,\ldots,0,1,0,\ldots,0$ as $x$ takes values
$1,2,\ldots,m$ modulo $m$ (with the 1 in the $a^{th}$ place). So we can
express $d(x)$ as a rational linear combination of the functions $k(x-a)$ for
$0\leq a<m$. This implies that we can express $d(x)$ as a rational linear
combination of functions of the form $\gcd(x-n_{i},m_{i})$ where $m_{i}%
=\frac{m}{d_{T}}$ for some $T\subset S$. If $m_{i}=1$ then we can replace
$\gcd(x-n_{i},m_{i})$ by the constant 1. Also, we can assume $0\leq
n_{i}<m_{i}$. Finally, using the fact that%
\[
\sum_{a=0}^{m_{i}-1}\gcd(x-a,m_{i})
\]
is a constant function, we can assume that $0<n_{i}<m_{i}$, provided we add a
constant term into our expression for $d(x)$. This completes the proof of
Theorem 2. Note that the proof shows that we can assume that the denominators
of the rational coefficients which appear in the expression for $d(x)$ divide
the constant $k(m)$.

\end{document}